\documentclass[a4paper,12pt,leqno]{amsart}
\usepackage{amsmath,amssymb,mathrsfs}
\usepackage[matrix,arrow]{xy} 
\usepackage{xcolor,hyperref}
\hypersetup{colorlinks=true,citecolor=black,linkcolor=black}
\newcommand{\h}{\hbox}
\newcommand{\q}{\quad}

\newcommand{\ms}{\par\medskip}
\newcommand{\sk}{\par\smallskip}

\newcommand{\msn}{\par\medskip\noindent}

\newcommand{\ges}{\geqslant}
\newcommand{\les}{\leqslant}
\newcommand{\one}{\hskip1pt}

\newcommand{\msum}{\hbox{$\sum$}}

\newcommand{\I}{{\mathcal I}}

\newcommand{\Hc}{{\mathcal H}}

\newcommand{\OO}{{\mathcal O}}
\newcommand{\Sc}{{\mathcal S}}
\newcommand{\X}{{\mathcal X}}

\newcommand{\Q}{{\mathbb Q}}
\newcommand{\C}{{\mathbb C}}
\newcommand{\N}{{\mathbb N}}
\newcommand{\R}{{\mathbb R}}
\newcommand{\RR}{{\mathbf R}}
\newcommand{\DD}{{\mathbf D}}
\newcommand{\Z}{{\mathbb Z}}

\newcommand{\Gr}{{\rm Gr}}
\newcommand{\al}{\alpha}
\newcommand{\be}{\beta}
\newcommand{\ga}{\gamma}
\newcommand{\Ga}{\Gamma}
\newcommand{\de}{\delta}
\newcommand{\De}{\Delta}
\newcommand{\la}{\lambda}
\newcommand{\La}{\Lambda}
\newcommand{\si}{\sigma}

\newcommand{\ep}{\varepsilon}
\newcommand{\qh}{\widehat{q}}
\newcommand{\qt}{\widetilde{q}}
\newcommand{\rt}{\widetilde{r}}

\newcommand{\sit}{\widetilde{\sigma}}
\newcommand{\Dt}{\widetilde{D}}

\newcommand{\Xo}{\overline{\mathcal X}}
\newcommand{\Yt}{\widetilde{Y}}
\newcommand{\Y}{\overline{Y}}
\newcommand{\Zo}{{}\,\overline{\!Z}{}}
\newcommand{\Om}{\Omega}
\newcommand{\dd}{\partial}
\newcommand{\ddd}{{\rm d}}

\newcommand{\Ff}{F_{\!f}}
\newcommand{\Gp}{\Gamma_{\!+}}
\newcommand{\Gf}{\Gamma_{\!f}}

\newcommand{\Sp}{{\rm Sp}}
\newcommand{\SpW}{{}^W\!\one {\rm Sp}}

\newcommand{\tos}{\,{\to}\,}
\newcommand{\col}{\,{:}\,}
\newcommand{\eq}{\,{=}\,}
\newcommand{\defs}{\,{:=}\,}
\newcommand{\nes}{\,{\ne}\,}
\newcommand{\ins}{\,{\in}\,}
\newcommand{\sst}{\,{\subset}\,}
\newcommand{\stm}{\,{\setminus}\,}

\newcommand{\less}{\,{\les}\,}
\newcommand{\sgt}{\,{>}\,}
\newcommand{\slt}{\,{<}\,}
\newcommand{\pl}{\one {+}\one}
\newcommand{\mi}{\one {-}\one}
\newcommand{\bl}{\bigl}
\newcommand{\br}{\bigr}

\newcommand{\ssb}{\raise.15ex\h{${\scriptscriptstyle\bullet}$}}
\newcommand{\ssc}{\,\raise.15ex\h{${\scriptstyle\circ}$}\,}

\newcommand{\into}{\hookrightarrow}
\newcommand{\simto}{\,\,\rlap{\hskip1.5mm\raise1.4mm\hbox{$\sim$}}\hbox{$\longrightarrow$}\,\,}

\makeatletter
\renewcommand\section{\@startsection{section}{1}{0pt}{-3ex plus -1ex minus -.2ex}{2.3ex plus.2ex}{\centering\normalfont\bfseries}}
\makeatother
\theoremstyle{plain}
\newtheorem{thm}{Theorem}[section]
\newtheorem{cor}{Corollary}[section]

\newtheorem{ithm}{Theorem}

\theoremstyle{definition}
\newtheorem{rem}{Remark}[section]

\begin{document}
\title[Spectral pairs of Newton non-degenerate polynomials]
{Spectral pairs of Newton non-degenerate polynomials with isolated singularities}
\author[S.-J. Jung]{Seung-Jo Jung}
\address{S.-J. Jung : Department of Mathematics Education, and Institute of Pure and Applied Mathematics, Jeonbuk National University, Jeonju, 54896, Korea}
\email{seungjo@jbnu.ac.kr}
\author[I.-K. Kim]{In-Kyun Kim}
\address{I.-K. Kim : June E Huh Center for Mathematical Challenges, Korea Institute for Advanced Study, 85 Hoegiro Dongdaemun-gu, Seoul 02455, Korea}
\email{soulcraw@kias.re.kr}
\author[M. Saito]{Morihiko Saito}
\address{M. Saito : RIMS Kyoto University, Kyoto 606-8502 Japan}
\email{msaito@kurims.kyoto-u.ac.jp}
\author[Y. Yoon]{Youngho Yoon}
\address{Y. Yoon : Department of Mathematics, Chungbuk National University, Cheongju-si, Chungcheongbuk-do, 28644, Korea}
\email{mathyyoon@gmail.com}
\thanks{This work was partially supported by National Research Foundation of Korea (the first author: RS-2026-25472637, the second author: NRF-2023R1A2C1003390 and NRF-2022M3C1C8094326, and the fourth author: RS-2023-00245670).}
\begin{abstract} We give a formula for the spectral pairs of Newton non-degenerate polynomials with isolated singularities. This is a revised version of a conjecture proposed by Steenbrink fifty years ago. Although this conjecture fails in the non-simplicial Newton polytope case, we can find the correction term using the intersection complexes and the primitive intersection cohomologies of toric varieties associated with non-simplicial compact faces. To verify the correctness of the formula, we provide a sample code for Singular which computes the spectral pairs and compares the output with that of sppairs in Singular for the three variable case.
\end{abstract}
\maketitle

\section*{Introduction}

Let $f\in\C\{x\}$ be a convergent power series having an isolated singularity at 0 with $f(0)\eq0$, where $x\eq(x_1,\dots,x_n)$ is a local coordinate system. We may assume $f$ is a polynomial by the finite determinacy, see for instance \cite{GLS}. We have the Milnor fiber $\Ff$, whose reduced cohomology is called the vanishing cohomology, and has the canonical mixed Hodge structure with the action of the semi-simple part of the monodromy. The {\it spectrum\one} $\Sp_f(t)$ is a fractional power polynomial defined by using the dimensions of graded pieces of the Hodge filtration $F$ on eigenspaces of the semi-simple part of the monodromy. Combining this construction with the information of {\it weights,} we can define the {\it weighted spectrum}
\begin{equation*}
\SpW_f(t,u)\eq\msum_{i=1}^{\mu_f}\,t^{\al_i}u^{w_i},
\end{equation*}
as the generating polynomial of {\it spectral pairs\one} $\{(\al_i,w_i)\}_{i\in[1,\one\mu_f]}$, see \cite{St1} and also Section~\ref{S1}. These are quite useful invariants, and have been studied by many people.
\sk
There is an explicit formula for the spectrum in the {\it convenient Newton non-degenerate\one} case, see \cite{St1}, \cite{exp}. Indeed, let $C(\si)\sst\R^n$ be the cone generated by a compact face $\si$ of the Newton polytope $\Gp(f)$, see \cite{Ko} and also \cite{JKSY1}. Let $p_{\si}(t)$ be the Hilbert-Poincar\'e series of the graded vector space $A_{\si}$ spanned by the monomials $x^{\nu}$ for $\nu\ins C(\si)\cap\N^n$ with respect to the Newton filtration $V'_N$, see \cite[Section 2]{simp}. We write $\tau\less\si$ if $\tau$ is a face of $\si$. Let $\Gf$ be the union of the compact faces of $\Gp(f)$. Set
\begin{equation} \label{1}
\begin{aligned}
d(\si):={}&\dim C(\si),\q c(\si)\defs n\mi d(\si)\q(\si\less\Gf),\\ \qh_{\si}(t):={}&(1-t)^{d(\si)}p_{\si}(t),\\ q_{\si}(t):={}&\msum_{\emptyset\les\tau\les\si}\,(-1)^{d(\si)-d(\tau)}\qh_{\tau}(t)\q\q\h{so that}\\ \qh_{\si}(t)\eq{}&\msum_{\emptyset\les\tau\les\si}\,q_{\tau}(t),
\end{aligned}
\end{equation}
where the summations are taken over all faces $\tau$ of $\si$ including $\si$ and $\emptyset$, see Remark~\ref{R4.1}.
\sk
We have the Steenbrink formula for the Poincar\'e polynomial $p_{B_{\Ga}}(t)$ of the graded quotients $B_{\Ga}\defs\Gr_{V_N}^{\ssb}\bl(\C\{x\}/(\dd f)\br)$ of the Newton filtration $V_N$ on the Jacobian ring $\C\{x\}/(\dd f)$\,:
\begin{equation} \label{2}
\begin{aligned}
p_{B_{\Ga}}(t)={}&\msum_{\emptyset\les\si\les\Gf}\,(-1)^{n-d(\si)}\,(1{-}t)^{k(\si)-d(\si)}\one\qh_{\si}(t)\\={}&\msum_{\emptyset\les\si\les\Gf}\,r_{\!\si}(t)q_{\si}(t)\q\q\q\q\q\q\,\,\h{with}\\ r_{\!\si}(t):={}&\msum_{\si\les\tau\les\Gf}\,(-1)^{n-k(\tau)}\,(t{-}1)^{k(\tau)-d(\tau)},
\end{aligned}
\end{equation}
where $q_{\emptyset}(t)\eq1$, and $r_{\!\si}(t)$ is called the {\it combinatorial polynomial,} see \cite[Theorem\,\,5.7]{St1} and also \cite[2.1]{simp}. This holds {\it without\one} assuming the primeness of compact faces. Note that \eqref{2} implies a formula for the spectrum using \cite{exp} and \cite[Remark 2.1f]{simp}.
\sk
It is known that the above assertion cannot be extended immediately to the spectral pairs unless the compact faces of the Newton polytope $\Gp(f)$ are {\it simplicial,} see \cite{Da}. It is however recently discovered that Steenbrink's conjectural formula for spectral pairs in the non-simplicial case can be modified by using the primitive cohomology of toric varieties corresponding to non-simplicial compact faces. Assume for simplicity that all compact faces $\si\less\Gf$ are {\it prime\one} (that is, for any vertex $v\ins\si$, the number of maximal-dimensional proper faces of $\si$ containing $v$ is equal to $\dim\si$, or equivalently the dual fan of $\si$ is simplicial; here it is enough to assume the condition only for maximal-dimensional compact faces). We have the following.

\begin{ithm} \label{T1}
Assume $f$ is Newton non-degenerate, has an isolated singularity at $0$, and all compact faces $\si\less\Gf$ are prime. Set $\rho_{\si,k}\defs\dim H_{\rm prim}^{2k}(\Xo_{\si})$ with $\Xo_{\si}$ the projective toric variety corresponding to $\si\less\Gf$. Steenbrink's conjecture for spectral pairs holds if we modify $q_{\sigma}(t)$ by subtracting $\rho_{\si}(t)\defs\msum_{k\ges1}\,\rho_{\si,k}\one t^k$ from it for any non-simplicial face $\si\less\Gf$ and adding $\msum_{\si}\,r_{\!\si}(t)\rho_{\si}(t)$ to the combinatorial polynomial $r_{\emptyset}(t)$ of $\emptyset$ up to sign. More precisely, setting
\begin{equation} \label{3}
\begin{aligned}
\qt_{\si}(t):={}&q_{\si}(t)\mi(-1)^{d_{\si}}\rho_{\si}(t),\\ \rt_{\emptyset}(t):={}&r_{\emptyset}(t)\pl\msum_{\si\les\Gf}\,(-1)^{d_{\si}}r_{\!\si}(t)\rho_{\si}(t),
\end{aligned}
\end{equation}
with $\rho_{\emptyset}(t)\eq0$ and $\rt_{\!\si}(t):=r_{\!\si}(t)$ for $ \si\nes\emptyset$, we have
\begin{equation} \label{4}
\SpW_f(t,u)=\msum_{\emptyset\les\si\les\Gf}\,\rt_{\!\si}(tu^2)\one\qt_{\si}(t)\one u^{d(\si)-1}.
\end{equation}
\end{ithm}

For the proof we show Steenbrink's formula for spectrum by using the {\it descent theorem\one} of nearby cycle formula for Newton non-degenerate polynomials (see \cite[Theorem 1]{des}). We construct in this paper an equivariant version of Griffiths theorem \cite{Gr} for toric varieties, which is different from the one in \cite{DK}, see Sections~\ref{S2} and \ref{3}. Using this, we can show also the {\it symmetry\one} of combinatorial polynomials $r_{\!\si}(t)$ except for $\si\eq\emptyset$ in the prime compact face case. (To examine the formula, we calculate $q_{\si}(t)$ using a computer in a case where $\si$ is homeomorphic to $I^{n-1}$ with $I\sst\R$ a closed interval for $n\eq4$ or 5. This is a natural extension of ``baby example" in \cite{Da}.)
\sk
The non-prime case is much more complicated, see Section~\ref{S6} below. We give a code for Singular computing the spectral pairs, and compare the output with that of ``sppairs" in Singular, see Section~\ref{S7}.

\tableofcontents
\numberwithin{equation}{section}

\section{Spectral pairs} \label{S1}

With the notation and assumption in the introduction the {\it spectral pairs\one} consist of $\mu_f$ pairs of rational numbers $(\al_i,w_i)$ ($i\ins[1,\mu_f]$) with $\al_i$ the spectral numbers and $w_i$ the modified weights. In this paper, the weights are shifted by $-1$ when $\al_i\ins\Z$. So the center of symmetry of the weights is always $n{-}1$. More precisely, we define the spectral pairs in this paper as follows (using the symmetry of spectral numbers):
\begin{equation*}
\begin{aligned}
&\#\{\one i\mid(\al_i,w_i)\eq(\al,w)\}=\dim\Gr_F^p\Gr^W_{w+\de_{\la,1}}H^{n-1}(\Ff,\C)_{\la}\\ &\q\q\h{with}\q\q p\eq\lfloor\al\rfloor,\q\la\eq e^{2\pi\sqrt{-1}\one\al},
\end{aligned}
\end{equation*}
where $\de_{\la,1}\eq1$ if $\la\eq1$, and 0 otherwise.

\begin{rem} \label{R1.1}
The above definition coincides with the one in \cite{St1} essentially, and is different from the one in \cite{SSS} by the shift of spectral numbers by $1$ and also by the change of modified weights under the involution of $\Z$ defined by $w\mapsto 2n\mi2\mi w$, where the {\it symmetry of modified weights\one} with center $n{-}1$ is used. Here it is actually better to change the spectral numbers (instead of modified weights) applying the involution of $\Q$ defined by $\al\mapsto n\mi\al$. Here the {\it self-duality\one} of the mixed Hodge structure on the vanishing cohomology is used.
\end{rem}

\section{A toric version of Griffiths theorem} \label{S2}

Let $\Gf$ be the union of the compact faces of $\Gp(f)$. For $\si\less\Gf$, the number of 1-dimensional faces of $\si$ containing $v$, denoted by $n_{v,\si}$, is equal to $d_{\si}\defs\dim\si$ for any vertex $v$ of $\si$ if and only if $\si$ is prime. Note that the maximal-dimensional compact faces $\tau$ of $\Gp(f)$ correspond to the irreducible components $\Xo_{\tau}$ of the exceptional divisor of the blowup associated with the dual fan of $\Gp(f)$, where the {\it combinatorics\one} of the compact faces of $\Gp(f)$ is the same as that of the intersections of the $\Xo_{\tau}$, see for instance \cite[\S2.1]{des}. Set $\Y\defs\Xo_{\si}$, the toric variety corresponding to $\si$, $Y:=\X_{\si}\sst\Y$ (the affine open torus), and $D\defs\Y\stm Y$, the complement divisor.

\begin{rem} \label{R2.1}
The logarithmic differential forms $\Om_{\Y}^p(\log D)$ are globally {\it free\one} sheaves {\it without\one} assuming $\si$ is prime, see \cite[\S2.9]{DK}. We define the logarithmic forms using the open direct image from the smooth locus of $\Y$. Its codimension is at least 2, since $\Y$ is normal. In the prime case, for each vertex $v\ins\si$, there are monomials $m_{v,i}\defs x^{\nu^{(v,i)}}$ ($i\ins[1,d_{\si}]$) with $\nu^{(v,i)}\ins\Z^n$ linearly independent and such that the $v{+}\be\nu^{(v,i)}$ are contained in the 1-dimensional faces of $\si$ containing $v$ for any $\be\ins[0,\ep]$ (with $0\slt\ep\,{\ll}\,1$). The logarithmic 1-forms $\ddd m_{v,i}/m_{v,i}$ then give local free generators on the torus-invariant affine open subvariety containing the point corresponding to $v$, see \cite[Remark 2.1a]{des}. For another vertex $v'\ins\si$, we have similarly vectors $\nu^{(v',i)}\ins\Z^n$ ($i\ins[1,d_{\si}]$). These two sets of vectors are related by an invertible matrix with $\Q$-coefficients (since they span the $\R$-vector space of dimension $d_{\si}$ spanned by the differences of two points of $\si$ considering the minimal affine space containing $\si$ if necessary), and so are the corresponding local free generators around $v$ and $v'$ (by the same matrix). Here vertices of $\si$ are identified with the corresponding points of $\Y$ as above. So the global freeness follows in the prime case. In general, we consider the direct image of the logarithmic forms on the ``prime cutting" (see Section~\ref{S3}), which is globally free. So the freeness follows, since $\Y$ is normal (using the Hartogs-type lemma).
\end{rem}

\begin{rem} \label{R2.2}
The above freeness implies that toric varieties have only {\it rational singularities,} using a desingularization $\pi\col\Yt\tos\Y$ via a subdivision of the dual fan of $\si$ and the Grauert-Riemenschneider vanishing theorem, since $\pi_*\I_{\Dt}=\I_D$. Here $\I_D\sst\OO_{\Y}$ is the reduced ideal sheaf of $D\sst\Y$ consisting of $g\ins\OO_{\Y}$ such that $g|_D\eq0$ as a function on $D$, and similarly for $\I_{\Dt}\sst\OO_{\Yt}$ with $\Dt\defs\pi^{-1}D\sst\Yt$. Indeed, the dualizing sheaf of $\Y$ is isomorphic to the tensor product of $\I_D$ with the logarithmic forms of the highest degree, using the freeness of logarithmic forms, the normality of $\Yt$, and the open direct image from the smooth locus of $\Y$ together with the Hartogs-type lemma, where we can employ also the Cohen-Macaulay property of $\Y$, see for instance \cite{Da0}.
\end{rem}

The freeness of logarithmic forms implies a version of Griffiths theorem \cite{Gr} (which is quite different from \cite{DK}) as follows. Let $\Zo\sst\Y$ be the closure of the hypersurface $Z\sst Y$ defined by $f_{\si}$. (Note that $f_{\si}$ is non-degenerate by assumption). Set $U\defs Y\stm Z$, and $l(k\si)\defs|\Z^n\,{\cap}\,k\si|$, the number of lattice points contained in $k\si$ ($k\ins\N$). Note that the complement $U$ is affine, and its complement in $\Y$, which is $D\,{\cup}\,Z$, is a divisor with $V$-normal crossings. Put
\begin{equation} \label{2.1}
\begin{aligned}
\chi_F^p(U)&\defs\msum_{j\in\Z}\,(-1)^j\dim_{\C}\Gr_F^pH^j(U)\,\,\,(\h{similarly for}\,\,Y,Z),\\ \chi_F^p(Y,Z)&\defs\msum_{j\in\Z}\,(-1)^j\dim_{\C}\Gr_F^pH^j(Y,Z),\\ P_{\!\si}(t)&\defs\msum_{k\ge0}\,l(k\si)\one t^k,\q P^{>0}_{\!\si}(t)\defs\msum_{k>0}\,l(k\si)\one t^k,
\end{aligned}
\end{equation}
so that $P_{\!\si}(t)\eq P^{>0}_{\!\si}(t)\pl 1$. We have the following.

\begin{thm} \label{T2.1}
There are equalities of polynomials
\begin{equation} \label{2.2}
\begin{aligned}
(-1)^{d_{\si}}\,\msum_{k\in\N}\,\chi_F^{d_{\si}-k}(Y,Z)\one t^k&=(1{-}t)^{d_{\si}+1}P_{\!\si}(t),\\ (-1)^{d_{\si}}\,\msum_{k\in\N}\,\chi_F^{d_{\si}+1-k}(U)\one t^k&=(1{-}t)^{d_{\si}+1}P^{>0}_{\!\si}(t).
\end{aligned}
\end{equation}
\end{thm}

\begin{proof}
There are long exact sequences of mixed Hodge structures
\begin{equation} \label{2.3}
\begin{aligned}
\tos H^j(Y,Z)\tos H^j(Y)\tos H^j(Z)\tos H^{j+1}(Y,Z)\tos,\\ \tos H^j(Z)\tos H^{j+2}(Y)(1)\tos H^{j+2}(U)(1)\tos H^{j+1}(Z)\tos.
\end{aligned}
\end{equation}
(These follow also from \cite{mhm}.) They imply that
\begin{equation} \label{2.4}
\begin{aligned}
\msum_k\,\chi_F^{d_{\si}-k}(Y)\one t^k&-\msum_k\,\chi_F^{d_{\si}-k}(Y,Z)\one t^k=\msum_k\,\chi_F^{d_{\si}-k}(Z)\one t^k,\\ \msum_k\,\chi_F^{d_{\si}+1-k}(Y)\one t^k&-\msum_k\,\chi_F^{d_{\si}+1-k}(U)\one t^k=\msum_k\,\chi_F^{d_{\si}-k}(Z)\one t^k.
\end{aligned}
\end{equation}
We have moreover
\begin{equation*}
\msum_k\,\chi_F^{d_{\si}-k}(Y)\one t^k\eq\msum_k\,\chi_F^{d_{\si}+1-k}(Y)\one t^{k-1}\eq(t{-}1)^{d_{\si}},
\end{equation*}
since $Y\cong(\C^*)^{d_{\si}}$. So the first and second equalities of \eqref{2.2} are equivalent to each other. We will show the second.
\sk
Assume $\si$ is {\it prime\one} for the moment. There is a bifiltered complex $\bl(\Om_{\Y}^{\ssb}(\log D)(*\Zo);F,W'\br)$, where $W'$ is induced by the weight filtration on $\Om_{\Y}^{\ssb}(\log D)$ via the tensor product with $\OO_{\Y}(*\Zo)$, and $F$ is defined by
\begin{equation} \label{2.5}
F^p\Om_{\Y}^j(\log D)(*\Zo)=\begin{cases}\Om_{\Y}^j(\log D)\bl((j{-}p{+}1)\Zo\br)&\h{if}\,\,\,p\less j,\\ \,0&\h{if}\,\,\,p\sgt j.
\end{cases}
\end{equation}
Using GAGA and a double complex structure around an intersection point of $\Zo$ and $D$, we verify that this filtered complex is filtered quasi-isomorphic to the logarithmic complex along $D\cup\Zo$ endowed with the filtration defined by the (stupid) truncation $\sigma^{\ges\ssb}$.
\sk
On the other hand, we have
\begin{equation} \label{2.6}
H^j\bl(\Y,\OO_{\Y}(\Zo)\br)=H^j\bl(\Y,\OO_{\Y}(\si)\br)=\begin{cases}L(\si)&\h{if}\,\,\,j\eq0,\\ \,0&\h{if}\,\,\,j\sgt0,\end{cases}
\end{equation}
see \cite{De}, \cite{Da0}, \cite[Proposition 2.8]{DK}. (This holds without assuming that $\si$ is prime.) Here $L(\si)$ is the $\C$-vector space generated by $x^{\nu}$ for $\nu\ins\Z^n\,{\cap}\,\si$, and $\OO_{\Y}(\si)$ is the invertible subsheaf of $\OO_{\Y}(*D)$ such that for a rational function on $\Y$, the order of zero or pole as a rational section of this line bundle along the divisor $\Xo_{\si'}\sst\Y$ for each $(d_{\si}{-}1)$-dimensional face $\si'\slt\si$ is given by using a defining linear function of $\si'\sst\si$ appropriately shifted. (Note that these shifts correspond to a cone-wise linear function on the dual fan of $\si$, which describes an {\it invertible sheaf,} see \cite{De}, \cite{Da0}. They can define an invertible sheaf, since intersection points of $(d_{\si}{-}1)$-dimensional faces are {\it lattice points,} which are vertices of $\si$, so that the coherent sheaf is generated locally by a {\it single\one} local section. This condition corresponds to that the restrictions of the corresponding cone-wise linear function to the 1-dimensional cones corresponding to toric divisors can be extended to a linear function on each maximal-dimensional cone, where the latter condition is satisfied by definition.) The first equality holds, since $\Zo$ is the {\it zero locus\one} of a ``general" section of the {\it invertible sheaf\one} $\OO_{\Y}(\si)$. Here ``general" means that the {\it non-degeneracy\one} condition is satisfied.
\sk
The equality \eqref{2.2} now follows from the freeness of the logarithmic complex in view of \eqref{2.5} and \eqref{2.6}. Indeed, concerning the multiplication by $(1{-}t)^{d_{\si}+1}$ on the right-hand side of \eqref{2.2}, the one by $(1{-}t)^{d_{\si}}$ corresponds to the passage from $\OO_{\Y}(*\Zo)$ to the {\it shifted\one} logarithmic de~Rham complex using the freeness, where $(-1)^{d_{\si}}$ on the left-hand side is needed to take this shift of the de~Rham complex into account. The remaining multiplication by $1{-}t$ corresponds to the passage from the filtration $F$ to the graded quotients $\Gr_F^{\ssb}$. (Concerning the leading term with $k\eq1$, note that $F^{d_{\si}}\Om_{\Y}^{d_{\si}}(\log D)(*\Zo)\eq\OO_{\Y}(\Zo)$.) So Theorem~\ref{T2.1} is proved in the prime case.
\sk
In the general case, we take a ``prime cutting" $\sit$ of $\si$ (see \cite[Section 3]{ict}), which gives a quasi-smooth subdivision of the dual fan of $\si$, see also \cite{Da0}. We define the line bundle $\OO_{\Yt}(\sit)$ on $\Yt\defs\X_{\si'}$ by using the ``limit" $\ep_{\si}\tos0$ for any $\si\slt\De$ in the notation of \cite[Section 3]{ict} (where $\si$, $\sit$ in our paper are denoted by $\De$, $\De'$ respectively) so that
\begin{equation} \label{2.7}
\OO_{\Yt}(\sit)=\pi^*\OO_{\Y}(\si),\q\OO_{\Y}(\si)=\pi_*\OO_{\Yt}(\sit),
\end{equation}
with $\pi\col\Yt\tos\Y$ a natural morphism, see also \cite{DK}. These isomorphisms imply that
\begin{equation} \label{2.8}
\RR\pi_*\OO_{\Yt}(\sit)=\OO_{\Y}(\si),\,\,\,\h{hence}\,\,\,R^k\pi_*\OO_{\Yt}(\sit)=0\,\,\,\h{if}\,\,\,k\sgt0,
\end{equation}
(using a desingularization of $\Yt$), since $\Yt$, $\Y$ are normal and have only rational singularities, see Remark~\ref{R2.2}. So the assertion follows from the freeness of logarithmic forms and \eqref{2.5} on $\Yt$ using \eqref{2.6} on $\Y$ together with \eqref{2.7} and \eqref{2.8}. (Here one cannot employ \eqref{2.6} {\it on\one} $\Yt$, since $\pi^*\Zo$ is {\it not\one} ample.) This finishes the proof of Theorem~\ref{T2.1}.
\end{proof}

\begin{cor} \label{C2.1}
If $\si$ is prime, we have an equality of polynomials
\begin{equation} \label{2.9}
(-1)^{d_{\si}}\msum_{k\in\N}\,\chi_F^k(\Y\stm\Zo)\one t^k=q_{\si}^{(1)}\!(t)+(-1)^{d_{\si}},
\end{equation}
where $q_{\si}^{(\la)}\!(t)\defs\msum_{k\ins\Z}\,q_{\si,\be+k}\one t^{\be+k}$ with $\la\eq e^{-2\pi\sqrt{-1}\be}$ if $q_{\si}(t)\eq\msum_{\al\in\Q}\,q_{\si,\al}t^{\al}$.
\end{cor}

\begin{proof}
Using the decomposition \cite[(5.3.3)]{MSS} (applied to the stratification $\Y\eq\bigsqcup_{\tau\les\si}X_{\tau}$), the duality isomorphism of mixed Hodge structures
\begin{equation} \label{2.11}
\DD H^k(X_{\tau},\Zo\cap X_{\tau})=H_c^{2d_{\tau}-k}\bl(X_{\tau},(j_{\tau})_*\Q_{h,X_{\tau}\setminus\Zo}\br)(d_{\tau}),
\end{equation}
and Theorem~\ref{T2.1}, we get the equalities
\begin{equation} \label{2.10}
\begin{aligned}
&(-1)^{d_{\si}}\msum_{k\in\N}\,\chi_F^k(\Y\stm\Zo)\one t^k\\&=(-1)^{d_{\si}}\msum_{\emptyset<\tau\les\si}\,\msum_{k\in\N}\,\chi_{c,F}^k\bl(X_{\tau},(j_{\tau})_*\Q_{h,X_{\tau}\setminus\Zo}\br)\one t^k\\&=(-1)^{d_{\si}}\msum_{\emptyset<\tau\les\si}\,\msum_{k\in\N}\,\chi_F^{d_{\tau}-k}(X_{\tau},\Zo\cap X_{\tau})\one t^k\\&=\msum_{k\in\N}\,q_{\si,k}\one t^k+(-1)^{d_{\si}},
\end{aligned}
\end{equation}
in view of the definition of $q_{\si}(t)$, see \eqref{1}. (Note that $\Y$ is quasi-smooth, since $\si$ is prime). Here $\DD$ denotes the dual functor of mixed Hodge structures, and the $j_{\tau}\col X_{\tau}\stm\Zo\into X_{\tau}$ are natural inclusions for $\tau\less\si$. We get $(-1)^{d_{\si}}\,\bl({=}\,(-1)^{d(\si)-1}\br)$ in the last term of \eqref{2.10}, since the summation in the definition of $q_{\si}(t)$ is taken over all faces of $\si$ including the {\it empty\one} one (see \eqref{1}), where $d(\emptyset)\eq0$ with $q_{\emptyset}(t)\eq\qh_{\emptyset}(t)\eq1$. Note also that $\OO_{\Y}(k\Zo)\cong\OO_{\Y}(k\si)$ for $k\ins\Z_{>0}$. This finishes the proof of Corollary~\ref{C2.1}.
\end{proof}

\section{An equivariant version} \label{S3}

We now study the {\it monodromy semisimple part equivariant\one} version of Section~\ref{S2}. Let $\si\less\Gf$ as in Section~\ref{S2}, and $\De$ be the convex hull of $\si\cup\{0\}$. We apply the ``prime cutting" as in Section~\ref{S2}, where the prime cutting at 0 is assumed to be given by a hyperplane parallel to the one containing $\si$. We get the toric variety $\Y_{\!\!\De'}$ associated with a prime cutting $\De'$ of $\De$. Here $\De'$ is infinitesimally close to $\De$, that is, we take essentially the ``limit" under $\ep_{\tau}\tos0$ for $\tau\slt\De$.
We have the isomorphism
\begin{equation*}
\Y_{\!\De}={\rm Proj}\,A_{\si}[t'],
\end{equation*}
so that ${\rm Spec}\,A_{\si}$ is the affine open subvariety containing 0 and its complement is $\Y\eq{\rm Proj}\,A_{\si}$. Here $\deg t'\eq\delta_{\si}^{-1}$ with $\delta_{\si}\in\N$ defined by
\begin{equation} \label{3.1}
\ell'_{\si}(\Z^n\cap V_{\si})=\tfrac{1}{\delta_{\si}}\Z,
\end{equation}
where $V_{\si}\sst\R^n$ is the $\R$-vector space spanned by $\si$ and $\ell'_{\si}$ is a linear function on $V_{\si}$ such that $\si\sst\ell_{\si}^{\prime-1}(1)$. The degrees of monomials contained in $A_{\si}$ are given by $\ell'_{\si}$. The monodromy semisimple part action, denoted by $T_s$, is given by multiplication by $e^{2\pi i/\delta_{\si}}$ on $t'$. This is compatible with the action on ${\rm Spec}\,A_{\si}$ defined by multiplication by $e^{-2\pi i\deg m}$ on monomials $m$ contained in $A_{\si}$. (It is the inverse of the Milnor or geometric monodromy, see also \cite{DS}.) This monodromy semisimple part action can be extended to $\Y_{\!\!\De'}$ so that its restriction to the exceptional divisor (and also to $\Y$) is the identity.
\sk
Let $\Zo_{\!\De}\sst\Y_{\!\De}$ be the closure of the affine hypersurface $Z_{\De}\sst Y_{\!\De}$ defined by $f_{\De}\defs f_{\si}\mi t'{}^{\delta_{\si}}$ in $\Y_{\!\De}={\rm Proj}\,A_{\si}[t']$. Here $Y_{\!\De}\sst\Y_{\!\De}$ is the affine open torus orbit isomorphic to $(\C^*)^{d_{\De}}$ and is contained in ${\rm Spec}\,A_{\si}\eq\Y_{\!\De}\stm\Y\eq\{t'\nes0\}$ with $\Zo_{\!\De}\cap{\rm Spec}\,A_{\si}$ defined by $f_{\si}\mi1$ (like a Milnor fiber of a homogeneous polynomial). The hypersurface $\Zo_{\!\De}$ is invariant by the action of the monodromy semisimple part $T_s$. For $\al\ins\Q$, set $k\defs\lceil\al\rceil$, $p\defs\lfloor\al\rfloor$, $\la\defs e^{2\pi \sqrt{-1}\al}$, and
\begin{equation} \label{3.2}
\begin{aligned}
\chi_{F,T_s}^{\al}(Z_{\De})&\defs\msum_{j\in\Z}\,(-1)^j\dim_{\C}\Gr_F^pH^j(Z_{\De})_{\la}\,\,\,(\h{similarly for}\,\,Y_{\De}),\\ \chi_{F,T_s}^{\al}(Y_{\De},Z_{\De})&\defs\msum_{j\in\Z}\,(-1)^j\dim_{\C}\Gr_F^pH^j(Y_{\De},Z_{\De})_{\la},\\ P_{\!\De,T_s}(t)&\defs\msum_{\al\ge0}\,P_{\!\De,T_s,\al}\one t^{\al}\q\h{with}\q P_{\!\De,T_s,\al}\defs l_{\overline{\la}}(k\De).
\end{aligned}
\end{equation}
Here $V_{\la}$ denotes the $\la$-eigenspace of $T_s$-action for any $\C$-vector space $V$ with $T_s$-action, and $l_{\overline{\la}}(k\De)$ is the dimension of the $\overline{\la}$-eigenspace of the $\C$-vector space spanned by the $x^{\nu}$ for $\nu\ins\Z^n\cap\De$ (and $\overline{\la}\eq\la^{-1}$ is the complex conjugate of $\la$). In the notation of \eqref{1}, we can verify the equalities
\begin{equation} \label{3.3}
(1{-}t)P_{\!\De,T_s}(t)=p_{\si}(t),\q (1{-}t)^{d_{\De}+1}P_{\!\De,T_s}(t)=\qh_{\si}(t),
\end{equation}
considering $k\De\stm(k{-}1)\De$ in $\R^n$ for the first equality.
We have the following.

\begin{thm} \label{T3.1}
There is a partial equality of fractional power polynomials
\begin{equation} \label{3.4}
(-1)^{d_{\si}}\,\msum_{\al\ge0}\,\chi_{F,T_s}^{d_{\De}-\al}(Z_{\De})\one t^{\al}=\qh_{\si}\q\h{\rm mod}\,\,\,\Z[t].
\end{equation}
\end{thm}

\begin{proof}
In the left-hand side of \eqref{3.4}, we can replace $\chi_{F,T_s}^{d_{\De}-\al}(Z_{\De})$ with $\chi_{F,T_s}^{d_{\De}-\al}(Y_{\De},Z_{\De})$ and $(-1)^{d_{\si}}$ with $(-1)^{d_{\De}}$ (where $d_{\si}\eq d_{\De}{-}1$) using the long exact sequence associated with the relative cohomology. (Indeed, the $T_s$-action on $H^{\ssb}(Y_{\De})$ is the identity considering a homotopy defined by multiplication by $e^{2\pi iu/\de_{\si}}$ on $t'$ for $u\ins[0,1]$, see the explanation of $T_s$-action after \eqref{3.1}.) 
Although the pullback of the invertible sheaf $\OO_{\Y_{\!\De}}(\De)$ by $\pi\col\Y_{\!\!\De'}\to\Y_{\!\De}$ is {\it not\one} ample so that Artin's vanishing theorem cannot be applied, we have
\begin{equation} \label{3.5}
R^j\pi_*\pi^*\OO_{\Y_{\!\De}}(\De)\eq0\q(j\sgt0),\q\pi_*\pi^*\OO_{\Y_{\!\De}}(\De)\eq\OO_{\Y_{\!\De}}(\De).
\end{equation}
by the same argument as in the proof of Theorem~\ref{T2.1}. Hence the assertion follows from the same argument as there, taking the action of $T_s$ into account. The equality holds only mod $\Z[t]$ because of the problem for the empty set $\emptyset\slt\si$. This finishes the proof of Theorem\,\,\ref{T3.1}.
\end{proof}

\begin{cor} \label{C3.1}
If $\si$ is prime, we have a partial equality of fractional power polynomials
\begin{equation} \label{3.6}
(-1)^{d_{\De}-1}\,\msum_{\al\ge0}\,\chi_{F,T_s}^{d_{\De}-\al}(\Zo_{\!\De})\one t^{\al}=q_{\si}(t)\q\h{\rm mod}\,\,\,\Z[t],
\end{equation}
together with a partial symmetry
\begin{equation} \label{3.7}
q_{\si}(t)=q_{\si}(t^{-1})\one t^{d_{\si}+1}\q{\rm mod}\,\,\,\,\Z[t].
\end{equation}
\end{cor}

\begin{proof}
Since $\si$ is prime, the varieties $\Zo_{\!\De}$, $\Zo$ are quasi-smooth, that is, these are $V$-manifolds (although $\Y_{\!\De}$ is not). Indeed, they are transversal to the stratifications of $\Y_{\!\De}\stm\{0\}$, $\Y$. (Note that $\Zo_{\!\De}\sst\Y_{\!\De}$ is locally defined by $h_{\si}{-}t''{}^{\de_{\si}}$ near a point of $\{t'\eq0\}\sst\Y_{\!\De}$, where $h_{\si}$, $t''$ are locally obtained by modifying $f_{\si}$, $t'$ appropriately.) So we have the duality isomorphisms
\begin{equation} \label{3.8}
\begin{aligned}
\DD H^k(\Zo_{\!\De})_{\la}&=H^{2d_{\si}-k}_c(\Zo_{\!\De})_{\overline{\la}}(d_{\si}),\\ \DD H^k(Z_{\De})_{\la}&=H^{2d_{\si}-k}_c(Z_{\De})_{\overline{\la}}(d_{\si}),
\end{aligned}
\end{equation}
where $\DD$ denotes the dual functor of mixed $\C$-Hodge structures. The second isomorphism holds with $\De$ replaced by its compact faces containing 0, except $\{0\}\slt\De$ corresponding to $\emptyset\slt\si$. The first assertion of Corollary~\ref{C3.1} then follows from Theorem\,\,\ref{T3.1} together with the last equality of \eqref{1} by using an argument similar to \cite[(5.3.3)]{MSS} (taking the monodromy action into account).
\sk
Since the variety $\Zo_{\!\De}$ is quasi-smooth (and $d_{\si}\eq d_{\De}{-}1$), the second assertion follows from the first using the first isomorphism of \eqref{3.8}. This finishes the proof of Corollary~\ref{C3.1}.
\end{proof}

\section{Symmetry of combinatorial polynomials} \label{S4}

In the {\it prime compact face\one} case, we have the {\it symmetry of combinatorial polynomials\one} $r_{\!\si}(t)$ for $\si\nes\emptyset $ (which is deduced from the one for the spectrum $\Sp_f(t)$) as follows.

\begin{thm} \label{T4.1}
Let $f$ be a Newton non-degenerate convenient polynomial whose Newton polytope has only prime compact faces. Then the symmetry of combinatorial polynomials $r_{\!\si}$ holds for $\si\nes\emptyset:$
\begin{equation} \label{4.1}
r_{\!\si}(t)\eq r_{\!\si}(t^{-1})\one t^{c(\si)}\q(\emptyset\slt\forall\,\si\less\Gf).
\end{equation}
\end{thm}

\begin{proof}
Substituting $x_i^{e_i}$ for $x_i$ in $f$ with $e_i$ sufficiently large integers, we may assume that the degree of each $f_{\si}$ is quite large, since this substitution does not change the combinatorial polynomials $r_{\!\si}(t)$. Writing $r_{\!\si}(t)\eq\msum_{i=0}^{c(\si)}r_{\!\si,i}t^i$ ($r_{\!\si,i}\ins\Z$), the assertion is equivalent to
\begin{equation} \label{4.2}
r_{\!\si,i}\eq r_{\!\si,c(\si)-i}\q\h{if}\,\,\,\,i\ins[0,c(\si)],\,\si\nes\emptyset.
\end{equation}
\sk
Let $\Sc_f$ be the set of non-empty compact faces of $\Gf$. Let $\La'\sst\La\defs\Sc_f{\times}\N$ be the subset consisting of $(\si,i)\ins\Sc_f{\times}\N$ with \eqref{4.2} unsatisfied. Assume $\La'\nes\emptyset$. Set
\begin{equation*}
G(t)\defs\msum_{(\si,i)\in\La'}\,r_{\!\si,i}\one t^iq_{\si}(t).
\end{equation*}
Since the $q_{\si}(t)$ has a symmetry mod $\Z[t]$ with central degree $d(\si)/2$ and $r_{\emptyset}(t)\ins\Z[t]$, we get the symmetry
\begin{equation} \label{4.3}
G(t)\eq G(t^{-1})\one t^n\q\h{mod}\,\,\,\Z[t],
\end{equation}
from the one of $\Sp_f(t)$, see Corollary\,\,\ref{C3.1}.
\sk
Let $m_{\si}\eq x^{\nu^{(\si)}}$ be the monomial corresponding to the lowest degree term $t^{\al_{\si}}$ of $q_{\si}(t)$, where $\nu^{(\si)}$ is contained in the interior of the cone of $\si$, using Remark~\ref{R4.1}. The degree $\al_{\si}$ can be obtained by applying $\ell'_{\si'}$ to $\nu^{(\si)}$, where $\ell'_{\si'}$ is the linear function whose value on $\si'$ is 1 with $\si'$ a maximal-dimensional compact face of $\Gf$ containing $\si$, see \eqref{3.1}.
\sk
Let $i_0$ be the minimal element of the image of $\La'$ by the second projection $\La\tos\N$. Set $\Sc'_{f,i_0}\defs\La'\cap\Sc_f{\times}\{i_0\}$, and $\al'\defs\min\{\al_{\si}\mid\si\ins\Sc'_{f,i_0}\}$. Then \eqref{4.3} implies that
\begin{equation} \label{4.4}
\#\{\si\ins\Sc'_{f,i_0}\mid\al_{\si}\eq\al'\}>1.
\end{equation}
However, this multiplicity becomes 1 if we modify the big integers $e_i$ appropriately (replacing the maximal-dimensional compact face $\si'$ containing $\si$ if necessary). This is a contradiction. So the assertion follows. This finishes the proof of Theorem\,\,\ref{T4.1}.
\end{proof}

\begin{rem} \label{R4.1}
For a compact face $\si\less\Gf$, let $\Z_{\si}$ be the constant sheaf on $\si$ with stalks $\Z$. We have the equality
\begin{equation} \label{4.5}
[(j_{\si})_!\Z_{\si^{\circ}}]=\msum_{\tau\les\si}\,(-1)^{d_{\si}-d_{\tau}}[\Z_{\tau}],
\end{equation}
in the Grothendieck group $K_{0,c}(\Z_{\si})$ of finitely generated $\Z_{\si}$-modules which are constructible for the natural stratification of $\si$. Here $\si^{\circ}$ denotes the interior of $\si$ with $j_{\si}\col\si^{\circ}\into\si$ the inclusion if $d_{\si}\nes0$, and $\si^{\circ}\eq\si$ if $d_{\si}\eq0$. The equality \eqref{4.5} follows from the identity
\begin{equation} \label{4.6}
[\Z_{\si}]=\msum_{\tau\les\si}\,[(j_{\tau})_!\Z_{\tau^{\circ}}],
\end{equation}
using isomorphisms
\begin{equation} \label{4.7}
\DD(j_{\tau})_!\Z_{\tau^{\circ}}=\Z_{\tau}[d_{\tau}]\q\h{for}\,\,\,\tau\less\si.
\end{equation}
Here an orientation of each $\tau\less\si$ is chosen, and $\DD$ denotes the dual functor, see \cite{Ve}.
\sk
Applying the restriction functor $i_{\si'^{\circ}}^*$ by the inclusion $i_{\si'^{\circ}}\col\si'^{\circ}\into \si$ for any {\it proper\one} face $\si'\slt\si$, we then get from \eqref{4.5} the vanishing
\begin{equation} \label{4.8}
\phi_{\si,\si'}\defs\msum_{\si'\les\tau\les\si}\,(-1)^{d_{\tau}}=0\q\q\h{for any}\,\,\,\si'\slt\si.
\end{equation}
This gives the last equality of \eqref{1}, see also \cite{des} (where the Euler characteristic is used).
Note that the class of a finite free $\Z_{\si'^{\circ}}$-module in the Grothendieck group depends only on its {\it rank.} It is possible to get an {\it acyclic complex\one} using an argument similar to the construction of ``realization" functor in \cite[3.1.9]{BBD}. We have a similar assertion with $\si$ replaced by its cone in $\R^n$.
\end{rem}

\section{Proof of Theorem~\ref{T1}} \label{S5}

The assertion for the {\it non-integral\one} exponent part follows from the calculations in Section~\ref{S3}. For the {\it integral\one} exponent part, we can apply Corollary~\ref{C2.1}, since $\Zo\sst\Y$ is quasi-smooth and is an ample divisor. Indeed, we have the long exact sequence
\begin{equation} \label{5.1}
\tos H^{j-2}(\Zo)(-1)\,{\buildrel{\!\ga}\over\to}\,H^j(\Y)\tos H^j(\Y\stm\Zo)\tos H^{j-1}(\Zo)(-1)\,{\buildrel{\!\ga}\over\to},
\end{equation}
using the local cohomology sequence. This implies a (non-canonical) short exact sequence
\begin{equation} \label{5.2}
0\tos H_{\rm prim}^{d_{\si}}(\Y)\tos H^{d_{\si}}(\Y\stm\Zo)\tos H_{\rm st.prim}^{d_{\si}-1}(\Zo)(-1)\tos0,
\end{equation}
together with (non-canonical) isomorphisms
\begin{equation} \label{5.3}
H^j(\Y\stm\Zo)\cong H_{\rm prim}^j(\Y)\q(j\in[0,d_{\si}{-}1]),
\end{equation}
using self-dualities of $H^{\ssb}(\Y)$, $H^{\ssb}(\Zo)$, where $H_{\rm st.prim}^j(\Zo)\defs{\rm Coker}\bl(i^*\col H^j(\Y)\tos H^j(\Zo)\br)$. The latter is non-canonically isomorphic to the kernel of the dual Gysin morphism $\ga$ using the self-duality and also semisimplicity of pure Hodge structures. The hard and weak Lefschetz theorems are closely related to the restriction and Gysin morphisms of $i\col\Zo\into\Y$, where the compositions of the two morphisms coincide with the action of the first Chern class of the ample line bundle. The cohomology $H^j(\Y)$ vanishes for $j$ odd, and has type $(j/2,j/2)$ for $j$ even (using the stratification by torus orbits and the Grothendieck group) since $\Y$ is quasi-smooth, that is, a $V$-manifold. So the correction term is given by the primitive cohomology of positive even degrees, and the symmetry of the modified $\qt_{\si}(t)$ follows from Corollaries~\ref{C3.1} and \ref{C2.1} together with \eqref{5.2} and \eqref{5.3} using the Grothendieck group. (Note that $(-1)^{d_{\si}}$ in \eqref{2.10} corresponds to $H^0(\Y)$.) We then get the symmetry of the modified $\rt_{\emptyset}(t)$ by Theorem~\ref{T4.1} and the symmetry of the spectrum. This finishes the proof of Theorem\,\,\ref{T1}.

\section{Non-prime case} \label{S6}

For $\emptyset\less\tau\slt\si\less\Gf$ with $\si\nes\emptyset$, set $m'_{\si,\tau}(t)\defs1\pl m_{\si,\tau}(t)$ and
\begin{equation} \label{6.1}
\begin{aligned}
m_{\si,\tau}(t)&:=\msum_{j>0}\, m_{\si,\tau,j}\one t^j\in\Z[t]\one t\q\q\h{with}\\ m_{\si,\tau,j}&:={\rm rank}\,\Hc^{2j-d_{\si}}{\rm IC}_{X_{\si}}\Q|_{X_{\tau}}\q\h{if}\,\,\,\tau\nes\emptyset,\\ m_{\si,\emptyset}(t)&:=\rho_{\si}(t),\\
\iota_{\si}(t)&:=\msum_{\emptyset\les\tau<\si}\,(-1)^{d_{\si}-d_{\tau}}m_{\si,\tau}(t)\one\qh_{\tau}(t),\\ \qt_{\si}(t)&:=q_{\si}(t)+\iota_{\si}(t).
\end{aligned}
\end{equation}
Here $q_{\emptyset}(t)\eq\qt_{\emptyset}(t)\eq1$ with $d_{\emptyset}\eq{-}1$ and the $\rho_{\si}(t)$ are defined by using the primitive part of intersection cohomology. Using \eqref{1}, we get by increasing induction on $d_{\si}$ that
\begin{equation} \label{6.2}
q_{\si}(t)=\qt_{\si}(t)+\msum_{\emptyset\les\tau<\si}\,n_{\si,\tau}(t)\one\qt_{\tau}(t)\q\h{with}\,\,\,\,n_{\si,\tau}(t)\in\Z[t]\one t.
\end{equation}
For $\emptyset\less\tau\less\Gf$, set
\begin{equation*}
\rt^{\one\prime}_{\si}(t)=r_{\si}(t)+\msum_{\Gf\ges\eta>\si}\,r_{\eta}(t)\one n_{\eta,\si}(t).
\end{equation*}
We then get
\begin{equation} \label{6.3}
\begin{aligned}
&\msum_{\emptyset\les\si\les\Gf}\,\rt^{\one\prime}_{\si}(t)\qt_{\si}(t)-\msum_{\emptyset\les\si\les\Gf}\,r_{\si}(t)q_{\si}(t)\\&=\msum_{\emptyset\les\si\les\Gf}\,n'_{\si}(t)\one\qt_{\si}(t)\q\h{with}\,\,\,\,n'_{\si}(t)\in\Z[t]\one t^2.
\end{aligned}
\end{equation}
Setting $\rt_{\si}(t)\defs\rt^{\one\prime}_{\si}(t)\mi n'_{\si}(t)$, we have
\begin{equation} \label{6.4}
\msum_{\emptyset\les\si\les\Gf}\,\rt_{\si}(t)\qt_{\si}(t)=\msum_{\emptyset\les\si\les\Gf}\,r_{\si}(t)q_{\si}(t).
\end{equation}
This implies the following.

\begin{thm} \label{T6.1}
With the above notation, the equality~\eqref{4} is extended to the non-prime case. 
\end{thm}

\begin{proof}
The fractional polynomials $\qt_{\si}(t)$ are the spectrum of {\it pure\one} Hodge structures with $T_s$-action so that the symmetry holds. This can be verified as in the proof of Theorem~\ref{T1} using the equalities
\begin{equation} \label{6.5}
\begin{aligned}
&(-1)^{d_{\si}}\msum_{k\in\N}\,\chi_{{\rm IH},F}^k(\Y\stm\Zo)\one t^k\\&=(-1)^{d_{\si}}\msum_{\emptyset<\tau\les\si}\,\msum_{k\in\N}\,\chi_{c,F}^k\bl(X_{\tau},(j_{\tau})_*\Q_{h,X_{\tau}\setminus\Zo}\br)\one m'_{\si,\tau}(t)\one t^k\\&=(-1)^{d_{\si}}\msum_{\emptyset<\tau\les\si}\,\msum_{k\in\N}\,\chi_F^{d_{\tau}-k}(X_{\tau},\Zo\cap X_{\tau})\one m'_{\si,\tau}(t)\one t^k\\&=\msum_{\emptyset<\tau\les\si}\,(-1)^{d_{\si}-d_{\tau}}\qh_{\tau}^{(1)}\!(t)\one m'_{\si,\tau}(t)\\&=q_{\si}^{(1)}\!(t)+(-1)^{d_{\si}}+\msum_{\emptyset<\tau\les\si}\,(-1)^{d_{\si}-d_{\tau}}\qh_{\tau}^{(1)}\!(t)\one m_{\si,\tau}(t),
\end{aligned}
\end{equation}

\begin{equation} \label{6.6}
\begin{aligned}
&(-1)^{d_{\si}}\msum_{k\in\N}\,\chi_{{\rm IH},F,T_s}^{\be+k}(\Zo_{\!\De})\one t^{\be+k}\\&=(-1)^{d_{\si}}\msum_{\emptyset<\tau\les\si}\,\msum_{k\in\N}\,\chi_{c,F,T_s}^{\be+k}(\Zo_{\!\De}\cap X_{\tau,\De})\one m'_{\si,\tau}(t)\one t^{\be+k}\\&=(-1)^{d_{\si}}\msum_{\emptyset<\tau\les\si}\,\msum_{k\in\N}\,\chi_{F,T_s}^{d_{\tau}+1-\be-k}(\Zo_{\!\De}\cap X_{\tau,\De})\one m'_{\si,\tau}(t)\one t^{\be+k}\\&=\msum_{\emptyset<\tau\les\si}\,(-1)^{d_{\si}-d_{\tau}}\qh_{\tau}^{(\la)}\!(t)\one m'_{\si,\tau}(t)\\&=q_{\si}^{(\la)}\!(t)+\msum_{\emptyset<\tau\les\si}\,(-1)^{d_{\si}-d_{\tau}}\qh_{\tau}^{(\la)}\!(t)\one m_{\si,\tau}(t)\q{\rm mod}\,\,\,\Z[t].
\end{aligned}
\end{equation}
Here $\chi_{{\rm IH},F}^k$, $\chi_{{\rm IH},F,T_s}^k$ are defined by employing the intersection cohomologies instead of the usual ones, and $\qh_{\tau}^{(\la)}\!(t)$ is defined like $q_{\si}^{(\la)}\!(t)$ in Corollary~\ref{C2.1}. We use also the direct image of the constant sheaf under the finite \'etale morphism $Z_{\De}\stm Y\tos Y\stm Z$, which gives the local system $L_{\si}$ in \cite[Theorem 1]{des}. The argument in Section~\ref{S4} is not affected, since we have $m_{\si,\tau}(t),\rho_{\si}(t)$ in $\Z[t]\one t$ so that the modification of $q_{\si}(t)$ by $\iota_{\si}(t)$ does not change the minimal degree, assuming the $e_i$ in Section~\ref{S4} are sufficiently large.
\end{proof}

\section{Sample code} \label{S7}

We give a sample code of Singular as follows. This is still a very experimental one. The input vertices are written as in the second line. (Any three vertices must not lie on the same line, and every coordinate axis must contain a vertex.) We first determine the combinatorics.
\ms
\vbox{\fontsize{6pt}{3mm}\tt\pv@LIB "linalg.lib"; LIB "gmssing.lib"; ring R=0,(x,y,t),ds; int a,b,c,h,i,j,k,m,v,p,q,r,s,n,Ms,Mv,Mc,Mr,Mn,fl;@
\pv@intvec Inp=2,1,1, 1,2,1, 4,0,0, 0,4,0, 0,0,5;@
\pv@poly f,X,Z; number mm; list spp; matrix WS[3][5]; v=size(Inp) div 3; intmat In[v][3]; for(i=1; i<=v; i++)@
\pv@{for(j=1; j<=3; j++){In[i,j]=Inp[3*i+j-3]; if(In[i,j]>m){m=In[i,j];}}} Ms=3*m; Mv=Ms; Mr=Ms; Mc=4*m; Mn=Mc;@
\pv@intmat De[4][Mc]; matrix Co[3][1]; matrix Inv[3][3]; vector Id=[1,1,1]; intmat C[Mc][4]; matrix V[Ms][Mv];@
\pv@intmat ks[3][Mc]; matrix Da[5][2]; matrix S[Ms][4]; intmat SS[Ms][Mv]; matrix Dg[3][3]=1,0,0,0,1,0,0,0,1;@
\pv@intmat Cs[Mc][3]; intmat OCs[Mc][3]; matrix Qh[3][Mn]; intmat Inc[Ms][v]; matrix Q[3][Mn]; matrix Rs[3][Mc];@
\pv@matrix A[3][3]; matrix Y[v][1]; matrix L[2][Mc]; intmat vs[v][1]; matrix B[4][2]; sprintf("Input vertices:");@}
\sk
\vbox{\fontsize{6pt}{2.5mm}\tt\pv@for(i=1; i<=v; i++){sprintf(" V%s: (%s,%s,%s)",i,In[i,1],In[i,2],In[i,3]);} s=1; for(i=1; i<=v-2; i++)@
\pv@{for(j=i+1; j<=v-1; j++){for(k=j+1; k<=v; k++){ for(p=1;p<=3;p++){A[1,p]=In[i,p]; A[2,p]=In[j,p];@
\pv@A[3,p]=In[k,p];} if(det(A)!=0){ Co=inverse(A)*Id; Y=In*Co; for(p=1;p<=3;p++){Da[p,1]=Co[p,1];} n=0;@
\pv@for(q=1; q<=v && sum(Y,q)>=1;q++){ if(sum(Y,q)>1){Inc[s,q]=1;} else {Inc[s,q]=0; n++; vs[n,1]=q;}}@
\pv@if(q>v){fl=0; for(m=1; m<=s && fl==0; m++){for(p=1; p<=3 && Da[p,1]==S[m,p]; p++){;} if(p>3){fl=1;}}@
\pv@if(fl==0){for(p=1; p<=3; p++){S[s,p]=Co[p,1];} SS[s,1]=n; SS[s,n+2]=1;@
\pv@sprintf("S%s: 2-dim face orthogonal to (%s) and containing %s vertices:",s,Co,n); for(p=1; p<=n; p++)@
\pv@{SS[s,1+p]=vs[p,1]; sprintf(" V%s",vs[p,1]);} s++;}}}}}} for(i=1; i<=3; i++){for(p=1; p<=3; p++)@
\pv@{S[s,p]=Dg[i,p]; Co[p,1]=Dg[i,p];} n=0; for(q=1; q<=v; q++){if(In[q,i]==0){n++; vs[n,1]=q;}}@}
\sk
\vbox{\fontsize{6pt}{2.5mm}\tt\pv@sprintf("S%s: 2-dim face orthogonal to (%s) and containing %s vertices:",s,Co,n); for(p=1; p<=n; p++)@
\pv@{SS[s,1+p]=vs[p,1]; sprintf(" V%s",vs[p,1]);} s++;} s--; for(i=1; i<=s; i++){for(q=1; q<=v; q++)@
\pv@{X=0; for(p=1; p<=3; p++){X=X+S[i,p]*In[q,p];} if(i>s-3){X=X+1;} V[i,q]=X;}} c=0; for(i=1; i<s; i++)@
\pv@{for(j=i+1; j<=s; j++){n=0; for(q=1; q<=v; q++){if(V[i,q]==1 && V[j,q]==1){n++; vs[n,1]=q;}}@
\pv@if(n>=2){c++; C[c,1]=i; C[c,2]=j; C[c,3]=vs[1,1]; C[c,4]=vs[2,1]; sprintf(@
\pv@"C%s: 1-dim face contained in S%s, S%s, and containing V%s, V%s", c,i,j,vs[1,1],vs[2,1]);}@
\pv@if(n>2){sprintf("V%s, V%s, V%s on the same line!",vs[1,1],vs[2,1],vs[3,1]); exit;}}}@}
\msn
We then compute the polynomials $\qh_{\si}$, $q_{\si}$, $r_{\si}$ and the correction terms as below.
\ms
\vbox{\fontsize{6pt}{2.5mm}\tt\pv@for(a=1; a<=s-3; a++){S[a,4]=S[a,1]*x+S[a,2]*y+S[a,3]*t; Z=cleardenom(S[a,4]); De[1,a]=int(number(diff(Z,x))@
\pv@/number(S[a,1])); q=0; for(b=1; b<=c; b++){if(C[b,1]==a || C[b,2]==a){q++; Cs[q,1]=b; Cs[q,2]=C[b,3];@
\pv@Cs[q,3]=C[b,4];}} if(q!=SS[a,1]){printf("Vertex Number Error");} for(b=1; b<=q; b++){A[1,1]=x; A[1,2]=y;@
\pv@A[1,3]=t; for(i=2; i<=3; i++){for(p=1; p<=3; p++){A[i,p]=In[Cs[b,i],p];}} X=det(A); for(h=2; h<=SS[a,1]+@
\pv@1 && subst(X,x,In[SS[a,h],1],y,In[SS[a,h],2],t,In[SS[a,h],3])>=0; h++){;} if(h<=SS[a,1]+1){X=-X;}@
\pv@L[1,Cs[b,1]]=X; L[2,b]=X; for(h=2; h<=SS[a,1]+1 && subst(X,x,In[SS[a,h],1],y,In[SS[a,h],2],t,In[SS[a,h],3])>=@
\pv@0; h++){;} if(h<=SS[a,1]+1){printf("Sign Error!");}} for(i=0; i<=Mr; i++){for(j=0; j<=Mr; j++)@
\pv@{for(k=0; k<=Mr; k++){for(b=1; b<=q && subst(L[2,b],x,i,y,j,t,k)>=0; b++){;} if(b>q){m=int(number(@
\pv@subst(Z,x,i,y,j,t,k))); Qh[1,a]=Qh[1,a]+x^m;}}}} Qh[1,a]=jet(Qh[1,a]*(1-x^(De[1,a]))^3,3*De[1,a]);}@}
\sk
\vbox{\fontsize{6pt}{2.5mm}\tt\pv@for(b=1; b<=c; b++){X=L[1,b]; Co[1,1]=diff(X,x); Co[2,1]=diff(X,y); Co[3,1]=diff(X,t); for(p=3; Co[p,1]==0;@
\pv@p--){;} for(i=1; i<=2; i++){if(i<p){j=i;} else {j=i+1;} Da[1,i]=Co[j,1]/number(Co[p,1]); B[1,i]=@
\pv@In[C[b,3],j]; B[2,i]=In[C[b,4],j]; } for(i=1; i<=2; i++){B[i+2,1]=B[i,2]*x-B[i,1]*y;} if@
\pv@(subst(B[3,1],x,B[2,1],y,B[2,2])<0){B[3,1]=-B[3,1];} if(subst(B[4,1],x,B[1,1],y,B[1,2])<0){B[4,1]=-B[4,1];}@
\pv@X=((B[2,2]-B[1,2])*x+(B[1,1]-B[2,1])*y)/(B[1,1]*B[2,2]-B[1,2]*B[2,1]); Z=cleardenom(X); if(diff(Z,x)!=0)@
\pv@{De[2,b]=int(number(diff(Z,x))/number(diff(X,x)));} else {De[2,b]=int(number(diff(Z,y))/number(diff(X,y)));}@
\pv@if(De[2,b]<0){De[2,b]=-De[2,b]; Z=-Z;} for(i=0; i<=Mr; i++){for(j=0; j<=Mr; j++){mm=number(Da[1,1]*i+@
\pv@Da[1,2]*j); if(subst(B[3,1],x,i,y,j)>=0 && subst(B[4,1],x,i,y,j)>=0 && mm==int(mm)){m=int(number(@
\pv@subst(Z,x,i,y,j))); Qh[2,b]=Qh[2,b]+x^m;}}} Qh[2,b]=jet(Qh[2,b]*(1-x^(De[2,b]))^2,3*De[2,b]);}@}
\sk
\vbox{\fontsize{6pt}{2.5mm}\tt\pv@for(b=1; b<=v; b++){p=0; for(i=1; i<=3; i++){if(In[b,i]!=0){p++; vs[p,1]=i;}} if(p==3){De[3,b]=@
\pv@gcd(gcd(In[b,1],In[b,2]),In[b,3]);} if(p==2){De[3,b]=gcd(In[b,vs[1,1]],In[b,vs[2,1]]);} if(p==1)@
\pv@{De[3,b]=In[b,vs[1,1]];} for(i=0; i<De[3,b]; i++){Qh[3,b]=Qh[3,b]+x^i;}}@}
\sk
\vbox{\fontsize{6pt}{2.5mm}\tt\pv@for(a=1; a<=s-3; a++){Q[1,a]=Qh[1,a]; for(b=1; b<=c; b++){fl=0; for(i=1; i<=c; i++){if(C[b,1]==a@
\pv@|| C[b,2]==a){fl=1;}} if(fl==1){r=lcm(De[1,a],De[2,b]);@
\pv@Q[1,a]=subst(Q[1,a],x,x^(r div De[1,a]))-subst(Qh[2,b],x,x^(r div De[2,b])); De[1,a]=r;}}@
\pv@for(i=2; i<=SS[a,1]+1; i++){r=lcm(De[1,a],De[3,SS[a,i]]); Q[1,a]=subst(Q[1,a],x,x^(r div De[1,a]))+@
\pv@subst(Qh[3,SS[a,i]],x,x^(r div De[3,SS[a,i]])); De[1,a]=r;} Q[1,a]=Q[1,a]-1;}@
\pv@for(b=1; b<=c; b++){Q[2,b]=Qh[2,b]; for(i=3; i<=4; i++){r=lcm(De[2,b],De[3,C[b,i]]); Q[2,b]=@
\pv@subst(Q[2,b],x,x^(r div De[2,b]))-subst(Qh[3,C[b,i]],x,x^(r div De[3,C[b,i]])); De[2,b]=r;} Q[2,b]=@
\pv@Q[2,b]+1;} for(b=1; b<=v; b++){Q[3,b]=Qh[3,b]-1;}@}
\sk
\vbox{\fontsize{6pt}{2.5mm}\tt\pv@for(b=1; b<=c; b++){if(C[b,1]>s-3 || C[b,2]>s-3){ks[2,b]=2; Rs[2,b]=0;} else {ks[2,b]=3; Rs[2,b]=x+1;}}@
\pv@for(b=1; b<=v; b++){p=0; for(i=1; i<=3; i++){if(In[b,i]==0){p++;}} ks[3,b]=3-p; Rs[3,b]=(1-x)^(ks[3,b]@
\pv@-1); for(i=1; i<=s-3; i++){for(j=2; j<=SS[i,1]+1 && SS[i,j]!=b; j++){;} if(j<=SS[i,1]+1){Rs[3,b]=@
\pv@Rs[3,b]+1;}} for(i=1; i<=c; i++){if(C[i,3]==b || C[i,4]==b){Rs[3,b]=Rs[3,b]-(1-x)^(ks[2,i]-2);}}}@
\pv@X=s-4; for(i=1; i<=c; i++){X=X-(1-x)^(ks[2,i]-2);} for(i=1; i<=v; i++){X=X+(1-x)^(ks[3,i]-1);}@}
\sk
\vbox{\fontsize{6pt}{2.5mm}\tt\pv@for(a=1; a<=s-3; a++){Rs[1,a]=1; Z=subst(jet(Q[1,a],De[1,a])-jet(Q[1,a],De[1,a]-1),x,1)-subst(jet(@
\pv@Q[1,a],2*De[1,a])-jet(Q[1,a],2*De[1,a]-1),x,1); Q[1,a]=Q[1,a]-Z*x^(De[1,a]); X=X+Z*x;@
\pv@if(Z-SS[a,1]+3!=0){sprintf("Correction Error at S%s!!!!",a); exit;}}  for(i=1; i<=5; i++){De[4,i]=1;}@
\pv@for(a=1; a<=s-3; a++){r=lcm(De[4,3],De[1,a]); WS[1,3]=subst(WS[1,3],x,x^(r div De[4,3]))+subst@
\pv@(Q[1,a],x,x^(r div De[1,a])); De[4,3]=r;}@}
\msn
The combinatorial and modified $q$-polynomials and the weighted spectrum are as follows.
\ms
\vbox{\fontsize{6pt}{2.5mm}\tt\pv@for(b=1; b<=c; b++){if(Rs[2,b]!=0){p=int(number(diff(Rs[2,b],x))); r=lcm(De[4,4],De[2,b]);@
\pv@WS[1,4]=subst(WS[1,4],x,x^(r div De[4,4]))+p*x^r*subst(Q[2,b],x,x^(r div De[2,b])); De[4,4]=r;@
\pv@p=int(number(subst(Rs[2,b],x,0))); r=lcm(De[4,2],De[2,b]); WS[1,2]=subst(WS[1,2],x,x^(r div De[4,2]))+@
\pv@p*subst(Q[2,b],x,x^(r div De[2,b])); De[4,2]=r;}}@}
\sk
\vbox{\fontsize{6pt}{2.5mm}\tt\pv@for(b=1; b<=v; b++){if(Rs[3,b]!=0){p=int(number(diff(diff(Rs[3,b],x),x)/2)); r=lcm(De[4,5],De[3,b]);@
\pv@WS[1,5]=subst(WS[1,5],x,x^(r div De[4,5]))+p*x^(2*r)*subst(Q[3,b],x,x^(r div De[3,b])); De[4,5]=r; p=@
\pv@int(number(subst(diff(Rs[3,b],x),x,0))); r=lcm(De[4,3],De[3,b]); WS[1,3]=subst(WS[1,3],x,x^(r div De[4,3]))+@
\pv@p*x^r*subst(Q[3,b],x,x^(r div De[3,b])); De[4,3]=r; p=int(number(subst(Rs[3,b],x,0)));@
\pv@r=lcm(De[4,1],De[3,b]); WS[1,1]=subst(WS[1,1],x,x^(r div De[4,1]))+p*subst(Q[3,b],x,x^(r div De[3,b]));@
\pv@De[4,1]=r;}} p=int(number(diff(diff(X,x),x)/2)); WS[1,4]=WS[1,4]+p*x^(2*De[4,4]);@
\pv@p=int(number(subst(diff(X,x),x,0))); WS[1,2]=WS[1,2]+p*x^(De[4,2]);@}
\sk
\vbox{\fontsize{6pt}{2.5mm}\tt\pv@sprintf("Combinatorial and modified q-polynomials:");@
\pv@for(i=1; i<=s-3; i++){sprintf("S%s: %s, %s (x%s=t)",i,subst(Rs[1,i],x,t),Q[1,i],De[1,i]);}@
\pv@for(i=1; i<=c; i++){sprintf("C%s: %s, %s (x%s=t)",i,subst(Rs[2,i],x,t),Q[2,i],De[2,i]);}@
\pv@for(i=1; i<=v; i++){sprintf("V%s: %s, %s (x%s=t)",i,subst(Rs[3,i],x,t),Q[3,i],De[3,i]);}@
\pv@sprintf("O: %s, 1",subst(X,x,t)); sprintf("Weighted spectrum:");@
\pv@for(i=1; i<=5; i++){sprintf("wt %s: %s (x%s=t)",i-1,WS[1,i],De[4,i]);}@}
\msn
If one is interested in the comparison, it can be done as below (although it may take long).
\ms
\vbox{\fontsize{6pt}{2.5mm}\tt\pv@for(i=1; i<=v; i++){f=f+(-1)^(i+1)*((i+2) div 2)*x^(In[i,1])*y^(In[i,2])*t^(In[i,3]);} spp=sppairs(f);@
\pv@for(i=1; i<=size(spp[3]); i++){j=5-spp[2][i]; WS[2,j]=WS[2,j]+spp[3][i]*x^(int((spp[1][i]+1)*De[4,j]));};@
\pv@for(i=1; i<=5; i++){WS[3,i]=WS[1,i]-WS[2,i]; if(WS[3,i]!=0){printf("Comp Error!!!"); exit;}}@
\pv@printf("Comp Ok");@}
\msn
Note that every computation requires the previous ones. It is recommended however not to do the computation at one time by unifying the code, especially for the comparison. (It may be better to edit it as an input file for Singular asking a specialist of Singular if necessary.) If the comparison fails, it might be possible that one computed a {\it Newton degenerate\one} polynomial $f$. In this case it may be solved by changing the order of vertices, since the coefficients of polynomial depend on the order. It may happen for instance if the input is as follows:
{\small\tt\verb@1,0,1,@} {\small\tt\verb@0,1,1,@} {\small\tt\verb@2,1,0,@} {\small\tt\verb@4,0,0,@} {\small\tt\verb@0,4,0,@} {\small\tt\verb@1,2,0,@} {\small\tt\verb@0,0,3;@} (here one can see $f$ by typing {\small\tt\verb@f;@} before running {\small\tt\verb@spp=sppairs(f);@}). It is not recommended to input too large numbers or too many vertices so that some constants must be changed.

\end{document}